\documentclass[11pt,reqno]{amsart}

\newtheorem{theorem}{Theorem}[section]

\theoremstyle{definition}

\theoremstyle{remark}

\theoremstyle{remark}
\newtheorem{example}[theorem]{Example}

\makeatletter
\def\dashint{\operatorname%
{\,\,\text{\bf--}\kern-.98em\DOTSI\intop\ilimits@\!\!}}
\makeatother

\newcommand{\nliminf}{\operatornamewithlimits{\underline{lim}}}

\newcommand\cF{\mathcal{F}}

 \newcommand{\mysection}[1]{\section{#1}
 \setcounter{equation}{0}}

\title[Result of A. Novikov
and Girsanov's theorem]{A few comments on a result of A. Novikov
and Girsanov's theorem}
\author[]{N.V. Krylov} 
\address{127 Vincent Hall, University of Minnesota,
 Minneapolis, MN, 55455}
 
\email{nkrylov@umn.edu}
\date{}%


\begin{document}

\begin{abstract}
 
We give a simple proof that for a continuous 
local martingale $M_{t}$ 
$$
\nliminf_{\varepsilon\downarrow0}\varepsilon
\log Ee^{(1-\varepsilon)
\langle M\rangle_{\infty}/2}<\infty
\Longrightarrow E\exp(M_{\infty}-\langle
 M\rangle_{\infty}/2)=1.
$$
 
\end{abstract}

\maketitle

\mysection{Main Result}

  Let $(\Omega,\mathcal{F},P)$ be a complete 
probability space and let
$M_{t}$ be a continuous local martingale 
on $(\Omega,\mathcal{F},P)$, provided with an appropriate
filtration of sub $\sigma$-fields of $\cF$,
such that $\langle M\rangle=\langle M
\rangle_{\infty}<\infty$ (a.s.).
Define
$$
M=M_{\infty},\,\,\,\rho=\rho(M)=e^{M-
\langle M\rangle/2},\,\,\,\rho_{t}=\rho_{t}
(M)=e^{M_{t}-
\langle M\rangle_{t}/2}.
$$
The process $\rho_{t}$ is called an exponential martingale,
although it is not necessarily a martingale.
We will be discussing generalizations of the following 
celebrated
 result of A.~Novikov (1972), which gives a sufficient
condition for $\rho_{t}$  to be a martingale:
\begin{equation}
                                           \label{10.17.1}
Ee^{\langle M\rangle/2}<\infty\Longrightarrow 
E\rho=1.
\end{equation}
This result is quite important in many applications
related to absolute continuous change of probability measure
and, in particular, makes available Girsanov's theorem.

The original proof in \cite{No_72} is based
on knowing the distribution of the first
exit time of the Wiener process with constant drift
from a shifted positive half-axis.
 Some of other known proofs are even more
  involved (see, for instance, Section 8.1 in \cite{RY_99}).
The latest easier proofs and the history
   revolving around
Novikov's condition are found and thoroughly discussed in \cite{Ru_13}
and \cite{Ru_15}.
Here we   present a completely elementary
 proof of a result that is somewhat
stronger than \eqref{10.17.1}.
\begin{theorem}
                                           \label{theorem 10.16.1}
We have

\begin{equation}
                                           \label{12.9.1}
\nliminf_{\varepsilon\downarrow0}\varepsilon\log 
Ee^{(1-\varepsilon)
\langle M\rangle/2}<\infty\Longrightarrow E\rho=1.
\end{equation}
\end{theorem}

Proof. We start with two known facts (proved at the end
of the paper):
\begin{equation}
                              \label{12.9.2}
 E\rho\leq1;\,\,\,                                                                 
 \exists\varepsilon>0:Ee^{(1+\varepsilon)
\langle M\rangle/2}<\infty\Longrightarrow E\rho=1
\end{equation}
and show a solution of part of Problem 4.3.13
of \cite{Kr_95} following the hint to that problem.
Observe that for small enough $\varepsilon>0$
$$
Ee^{(1+\varepsilon)^{2}\langle (1-\varepsilon)M\rangle/2}
= Ee^{(1-\varepsilon^{2})^{2}\langle M\rangle/2}<\infty,
$$
which by (\ref{12.9.2})   implies that $E\rho((1-\varepsilon)M)=1$.

Then use H\"older's inequality and  
write that for small enough $\varepsilon>0$
and any constant $T\in(0,\infty)$
$$
1=E\rho((1-\varepsilon)M)=Ee^{(1-\varepsilon)
(M-\langle M\rangle/2)}
e^{(1-\varepsilon)\varepsilon\langle
 M\rangle/2}I_{\langle M\rangle\leq T}
$$
$$
+EI_{\langle M\rangle> T}e^{(1-\varepsilon)
(M-\langle M\rangle/2)}
e^{(1-\varepsilon)\varepsilon\langle M\rangle/2}
$$
$$
\leq(E\rho)^{1-\varepsilon}
(Ee^{(1-\varepsilon)\langle M\rangle/2}I_{\langle M\rangle\leq T})^
{\varepsilon}+
(E\rho I_{\langle M\rangle> T})^{1-\varepsilon}
(Ee^{(1-\varepsilon)\langle M\rangle/2})^
{\varepsilon}.
$$
 
As $\varepsilon\downarrow0$, we get $$1
\leq E\rho+{\rm const} \,
E\rho I_{\langle M\rangle> T},$$  which gives $1\leq E\rho$
after letting $T\rightarrow\infty$.  This
together with the first relation in (\ref{12.9.2})
 implies
our statement (\ref{12.9.1}). The theorem is proved.

  Assertion  (\ref{12.9.1})  
 is stronger than \eqref{10.17.1}. 
\begin{example}
                                       \label{example 10.17.1}
 Take a one-dimensional Wiener 
process $w_{t}$ and let $\tau$
be the first exit time of $w_{t}$
from $(-\pi,\pi)$. Take $\lambda\in(0,1/8)$ and observe that
$u_{\lambda}(x)=\cos(\sqrt{2\lambda}x)$
satisfies $(1/2)u''+\lambda u=0$
on $[-\pi,\pi]$.
It follows by It\^o's formula
that $m^{\lambda}_{t}:= u_{\lambda}(w_{t\wedge\tau})
\exp(\lambda(t\wedge\tau))$ is a martingale and since 
$u_{\lambda}$ is bounded away from zero ($\lambda<1/8$),
there is a $\varepsilon>0$ such that, for any $t<\infty$,
$\varepsilon E\exp(\lambda(t\wedge\tau))\leq
Em^{\lambda}_{t}= 1$. As $t\to\infty$, we get
$\varepsilon E\exp \lambda \tau \leq 1$,
which by the dominated convergence theorem
allows us to send $t\to\infty$ in $Em^{\lambda}_{t}= 1$
and obtain that
$$
\cos(\sqrt{2\lambda}\pi)E\exp (\lambda\tau)=1,\quad
E\exp (\lambda\tau)=\big[\cos(\sqrt{2\lambda}\pi)\big]^{-1}
$$
if $\lambda<1/8$ and then $E\exp (\lambda\tau)=\infty$
if $\lambda= 1/8$.
Therefore, for the martingale $M_{t}=w_{t\wedge\tau}/4$
one easily finds that as
$\varepsilon\downarrow0$
$$
[Ee^{(1-\varepsilon)\langle M\rangle/2)}
]^{\varepsilon}=
[Ee^{(1-\varepsilon)\tau/8}]^{\varepsilon}
=\big[\cos \frac{\sqrt{1-\varepsilon}}{2}\pi
\big]^{-\varepsilon}\rightarrow1,
$$
so that the assumption in (\ref{12.9.1}) is 
satisfied, whereas
$E\exp(\langle M\rangle/2)=\infty$ and  
 Novikov's criterion is not
applicable.
\end{example}

\mysection{Refined Novikov's conditions}

In \cite{No_79} Novikov (1979)  relaxes the  conditions
 in \eqref{10.17.1}
and shows that, for any constant $c\geq0$,
\begin{equation}
                                           \label{10.17.2}
Ee^{\langle M\rangle/2-c\langle M\rangle^{1/2}}
<\infty\Longrightarrow 
E\rho=1.
\end{equation}
This condition is applicable  in Example \ref{example 10.17.1},
although it is not very easy to see that. We need to know
the tail of the distribution of $\tau$. On the other hand,
the elementary inequality:
$c\langle M\rangle^{1/2}\leq (\varepsilon/2) 
\langle M\rangle+c^{2}\varepsilon ^{-1}$
implies that
$$
I:=Ee^{\langle M\rangle/2-c\langle M\rangle^{1/2}}
<\infty\Longrightarrow 
\nliminf_{\varepsilon\downarrow0}\varepsilon\log 
Ee^{(1-\varepsilon)
\langle M\rangle/2}\leq \log I-c^{2}<\infty.
$$
 
In the same article \cite{No_79} Novikov gives
a more elaborated condition
\begin{equation}
                                           \label{10.17.3}
Ee^{\langle M\rangle/2-g(\langle M\rangle^{1/2})} 
<\infty\Longrightarrow 
E\rho=1,
\end{equation}
where $g$ belongs to the lower Kolmogorov class,
for instance, 
$$
g(\tau)=\sqrt{2\tau\log\log\tau}
$$
for large $\tau$. This condition is, of course,
much weaker than ours, but checking it, generally,
 requires much more
knowledge about the distribution of $\langle M\rangle$.

On the one hand, \eqref{10.17.3} provides the most general
condition in terms of the distribution
of $\langle M\rangle$. But on the other hand, this
 distribution, generally,  has little to do with the equality
$E\rho=1$. Indeed, if $\tau$ is the first time the Wiener
process $w_{t}$ hits point 1 and $M_{t}=w_{t\wedge\tau}$,
then by Kazamaki's criterion (see \cite{Ka_77})
we have $E\rho=1$ and at the same time even $
E\sqrt{\langle M\rangle}=E\sqrt{\tau}=\infty$
not to mention any exponential moments.
By the way, everything said about \eqref{12.9.1}
 has its natural counterpart
for Kazamaki's criterion (see 
http://arxiv.org/abs/math/0207013).

The conditions described above usually are interesting
not exclusively in their own rights but in connection
with the problem of absolute continuity of
the distribution of a stochastic process with respect to the
Wiener measure. For instance, let $\xi_{t}$ be a (adapted) solution
of a stochastic equation $d\xi_{t}=dw_{t}
+b(\xi_{\cdot},t)\,dt
$, $\xi_{0}=0$, with  
   nonanticipating $b(\xi_{\cdot},t)$ and we are interested
to know when its distribution $\mu_{\xi_{\cdot}}$ on $C(0,1)$
is absolutely continuous with respect to the distribution
$\mu_{w _{\cdot}}$  on $C(0,1)$ of the Wiener process
$  w_{t}$.
  According to Theorem 6 of \cite{LSh_72}, $\mu_{\xi_{\cdot}}\ll
\mu_{w _{\cdot}}$ iff
\begin{equation}
                                               \label{2.1.2}
P\big(\int_{0}^{1}b^{2}(\xi_{\cdot},t)\,dt<\infty\big)=1.
\end{equation}
Under this condition for any nonnegative
measurable function $f(x_{\cdot})$ on $C(0,1)$
\begin{equation}
                                               \label{2.1.1}
Ef(\xi _{\cdot})=Ef(w^{x}_{\cdot})e^{\phi},
\quad \phi=\int_{0}^{1}b(w^{x}_{\cdot},t)\,dw_{t}
-(1/2)\int_{0}^{1}b^{2}(w^{x}_{\cdot},t)\,dt.
\end{equation}
In particular, $Ee^{\phi}=1$. Thus \eqref{2.1.2}
implies that $Ee^{\phi}=1$.
 More general
multidimensional
equations are treated in \cite{LSh_01}.  

In conclusion, for the sake of completeness
we prove \eqref{12.9.2} following 
the proof of Lemma 3 in \cite{LSh_72}.
This shows that unlike \cite{No_72} and \cite{No_79}
 no specific information
about the Wiener process is needed to get our main result.

That $E\rho\leq1$ follows from the fact that $\rho_{t}$
is a local martingale (It\^o's formula) and the fact that
$\rho_{t}\geq0$, so that it is a supermartingale,
has the limit as $t\to\infty$ and Fatou applies.

 Next,
for any $p,r >1$, stopping time $\tau$ such that $\rho_{t\wedge\tau}$
is a martingale, and  
$t\in[0,\infty)$, we easily obtain by the H{\"{o}}lder inequality 
$$ 
E\rho_{t\wedge\tau}^p \leq \bigg(E 
\rho_{t\wedge\tau}(prM)\bigg)^{1/r}
  \bigg(E\exp(\frac{rp-1}{2}p\frac{r}{r-1}
\langle M\rangle_{t\wedge\tau} )\bigg)^
{(r-1)/r}.
$$
By \thetag4, the first factor on the right is at most $1$.
The second factor can be made bounded uniformly with respect to 
$t $, since $\langle M\rangle_{t}\leq\langle M\rangle$
and for $p=1+\delta,\,r=1+\sqrt{\delta}$
and sufficiently small $\delta >0$, it is not difficult to see
that the coefficient of $\langle M\rangle_{t\wedge\tau}$ can be made smaller than
$(1+\varepsilon)/2$. Fix $\delta,p$, with these properties.
Then by Doob's inequality $E\sup_{t\leq\tau}\rho_{t}^{p}\leq N$,
where the constant $N$ is independent of $\tau$. This yields that
$E\sup_{t<\infty}\rho_{t}^{p}\leq N$,
the local martingale $\rho_{t}$ is bounded by a summable
function independent of $t$ and $E\rho=\lim_{n\to\infty}
E\rho_{\tau_{n}}=1$, where $\tau_{n}\to\infty$ is any
localizing sequence for $\rho_{t}$.

{\bf Acknowledgment}. The author is sincerely grateful
to the referee for the comments and suggestions
which certainly helped improve the presentation.

\end{document}